\documentclass[3p,10pt,a4paper,twoside,fleqn,sort&compress]{preprint}

\usepackage{amssymb,amsmath,latexsym}
\usepackage[varg]{pxfonts}
\usepackage{xcolor}        % <------
\usepackage{url}           % <------
\usepackage{graphicx}
\usepackage{float}
\usepackage{url}
\usepackage{lscape} %package for figure and table placement

\newtheorem{theorem}{Theorem}[section]

\newtheorem{definition}[theorem]{Definition}

\newtheorem{remark}[theorem]{Remark}

\newcommand{\FilPages}{1--16}

\newcommand{\FilReceived}{Received: 29 February 2024; Accepted: dd Month yyyy}

\newcommand{\NN}{{\mathbb N}}
\newcommand{\RR}{{\mathbb R}}

\def\nfrac#1#2{\mbox{\small $\displaystyle\frac{#1}{#2}$}}
\def\nnfrac#1#2{\mbox{\footnotesize $\displaystyle\frac{#1}{#2}$}}

\usepackage[T1]{fontenc}                         % ************ %
\usepackage{amsmath,amsfonts}                    % ************ %
\DeclareMathAlphabet{\mathbbmsl}{U}{bbm}{m}{sl}  % ************ %

\begin{document}

\title{The area of H\" ugelsch\" affer curves via Taylor series}

\author[SF]{Maja Petrovi\' c\corref{A1}}
\ead{majapet@sf.bg.ac.rs}
\author[ETF]{Branko Male\v sevi\' c}
\ead{malesevic@etf.rs}

\address[SF]{University of Belgrade, The Faculty of Transport and Traffic Engineering, Serbia}
\address[ETF]{University of Belgrade, School of Electrical Engineering, Serbia}
\newcommand{\AuthorNames}{Maja Petrovi\' c and Branko Male\v sevi\' c}

\newcommand{\FilMSC}{Primary 14H50; Secondary 33E05, 41A10}
\newcommand{\FilKeywords}{H\" ugelsch\" affer curves, Elliptic integral, Taylor approximations}
\cortext[A1]{* Correspoding author: Maja Petrovi\' c}

\begin{abstract}
In this paper, we give new Taylor approximative formulae for the area of the egg-shaped parts of
H\" ugelsch\" affer curves. Based on a parametrization of the H\" ugelsch\" affer curve, a formula for the area of
the egg-shaped part of such a curve is derived via elliptic integrals of the first and second kind. Furthermore, new
approximative formulae for calculating this area derived from standard and double Taylor approximations are given.
A representation of the value $\frac{1}{\pi}$ was also obtained using an appropriate series.
\end{abstract}

\maketitle
\makeatletter
\renewcommand\@makefnmark%
{\mbox{\textsuperscript{\normalfont\@thefnmark)}}}
\makeatother

\section{H\" ugelsch\" affer curve $\mathcal{F}_{\mbox{\scriptsize \boldmath $q$},\, egg}$}

\begin{figure}[!h]
  \centerline{
    \includegraphics*[width=0.99\textwidth]{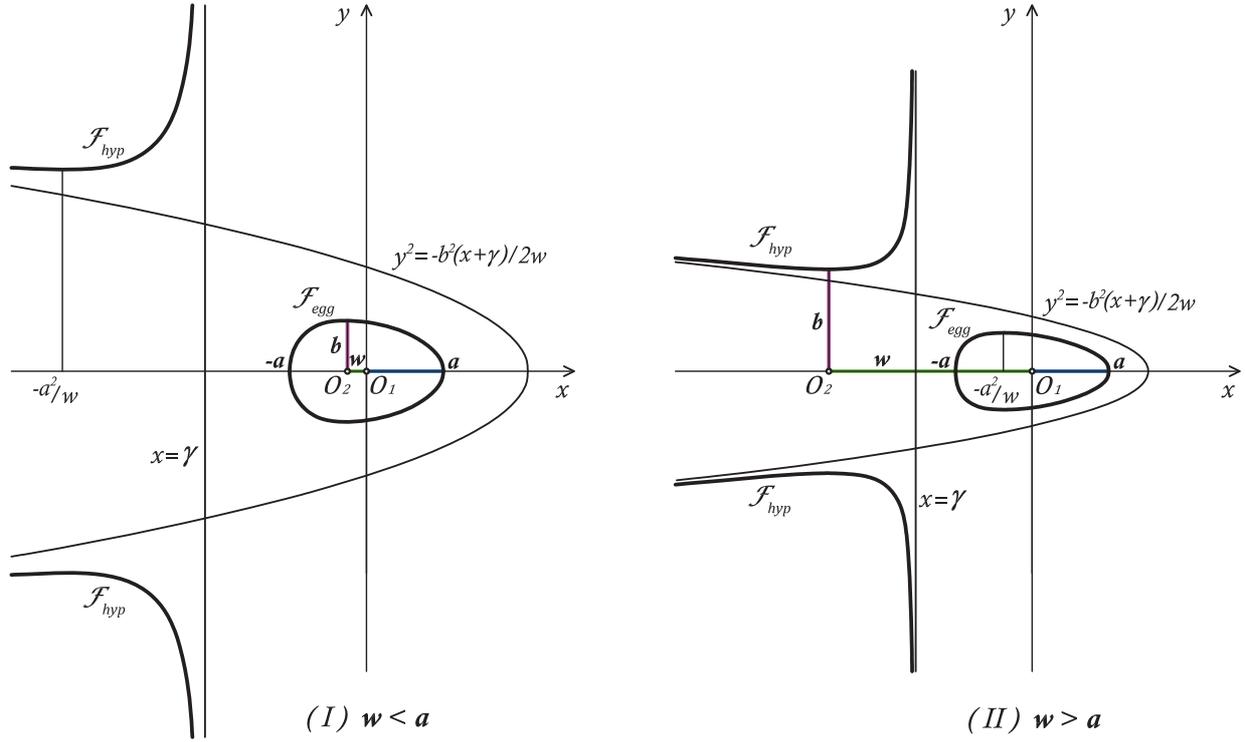}
  }
\caption{A H\" ugelsch\" affer curve $\mathcal{F} = \mathcal{F}_{egg} \cup \mathcal{F}_{hyp}$; \, Source: \copyright\,First Author}
\label{Fig.1}       % Give a unique label
\end{figure}

The H\" ugelsch\" affer curve \cite{FH 1944} is an algebraic cubic curve given by the following equation
\begin{equation}
\label{F}
\mathcal{F}: \quad 2wxy^2+b^2x^2+(a^2+w^2)y^2-a^2b^2=0,
\end{equation}
where $a, b, w > 0$. In the papers \cite{MP 2010}, \cite{moNG 2010a} and \cite{AADM 2023}, a decomposition of this cubic curve is described:
\begin{equation}
\label{UnionF}
\mathcal{F}
=
\mathcal{F}_{egg}
\cup
\mathcal{F}_{hyp}.
\end{equation}
The egg-shaped part $\mathcal{F}_{egg}$ of the curve is defined over $[-a,a]$, and the hyperbolic part $\mathcal{F}_{hyp}$ of the curve
(which consists of two branches) is defined over $(-\infty,\gamma)$, where $\gamma = - \frac{a^2+w^2}{2w}$, see Fig. \ref{Fig.1}.
It is easy to check that $\gamma <-a \Leftrightarrow (a-w)^2>0$.

%% FIGURE 1 %%%%%%%%%%%%%%%%%%%%%%%%%%%%%%%%%%%%%%%%%%%%%%%%%%%%%%%%%%%%%%%%%%%%%%%%%%%%%%%%%%%%%%%%

%\begin{equation}
%\label{ABC}
%\alpha=a,\;
%\beta=-a,\;
%\gamma=-\frac{a^2+w^2}{2w},\;
%H=\displaystyle\frac{b}{\sqrt{2w}}.
%\end{equation}
%as considered in \cite{MP 2010}.

\begin{figure}[htb]
\vspace{-0.3cc}
  \centerline{
    \includegraphics*[width=0.85\textwidth]{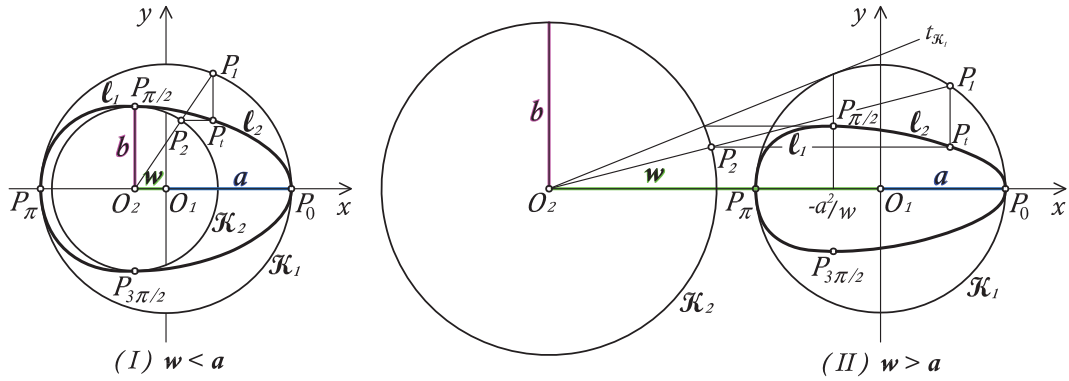}
  }
\caption{H\" ugelsch\" affer's  construction of an egg curve $\mathcal{F}_{\mbox{\scriptsize \boldmath $q$},\, egg}$; \, Source: \copyright\,First Author}
\label{Fig.2}       % Give a unique label
\end{figure}

\medskip
Let us consider just the non-degenerative cases of the H\" ugelsch\" affer cubic curve (\ref{F}) as in \cite{MP 2010} and \cite{AADM 2023}.
Let $w$ be the distance between the two circle centers when constructing the curve ($w=|O_1 O_2|$, see Fig.~\ref{Fig.2}).
Then, we consider the two cases $(I)$ $w < a$ and $(II)$ $w > a$. The abscissa \mbox{\boldmath $u$} at which the
H\" ugelsch\" affer curve reaches its extremes over the segment $[-a,a]$ is defined in \cite{AADM 2023} with the following formula:
\begin{equation}
\label{Uu}
\mbox{\boldmath $u$}
=
\left\{
\begin{array}{ccc}
-w     &:& w < a, \\[1.0 ex]
-a^2/w &:& w>a.
\end{array}
\right.
\end{equation}
Introducing the parameter \mbox{\boldmath $q$} as:
\begin{equation}
\label{QUu}
\mbox{\boldmath $q$}
=
\left\{
\begin{array}{ccc}
1     &:& w < a, \\[1.0 ex]
a/w   &:& w > a
\end{array}
\right\}
=
\frac{\sqrt{a}}{\sqrt{w}}\,
\frac{\min \{\! \sqrt{a}, \sqrt{w} \, \} }{\max \{\! \sqrt{a}, \sqrt{w} \, \}}\,.
\end{equation}
Then, it is true that $\mbox{\boldmath $q$} \in (0,1]$.

\smallskip
\noindent
The relationship between \mbox{\boldmath $u$} and \mbox{\boldmath $q$} (i.e. (\ref{Uu}) and (\ref{QUu})) is given as
\begin{equation}
\label{UuN}
\mbox{\boldmath $u$}
=
-\mbox{\boldmath $q$}^2w.
\end{equation}
Only case (I) of the cubic curve $\mathcal{F}_{egg}$ was considered by F. H\" ugelsch\" affer in the work \cite{FH 1944},
while case (II) was introduced by M. Petrovi\'c in the thesis \cite{MP 2010}.
% U ovom delu pokazujemo da je mogu\' ce ova dva slu\v caja posmatrati objedinjeno sa krivom
% \begin{equation}
% \label{Fq}
% \mathcal{F}_{\mbox{\scriptsize \boldmath $t$}}: \quad 2\mbox{\boldmath $t$}^2wxy^2+\mbox{\boldmath $t$}^2b^2x^2
% +(a^2+\mbox{\boldmath $t$}^4w^2)y^2-a^2b^2\mbox{\boldmath $t$}^2=0,
% \end{equation}
% sa odgovaraju\' com vredno\v s\' cu realnog parametra $\mbox{\boldmath $t$}$ za koju je
% \begin{equation}
% \label{Eq}
% \mathcal{F} = \mathcal{F}_{\mbox{\scriptsize \boldmath $t$}}.
% \end{equation}
% Uslov (\ref{Eq}) je ispunjen ako i samo ako
% $$
% \frac{2 \mbox{\boldmath $t$}^2 w}{2 w}
% =
% \frac{\mbox{\boldmath $t$}^2 b^2}{b^2}
% =
% \frac{a^2 + \mbox{\boldmath $t$}^4w^2}{a^2+w^2}
% =
% \frac{a^2b^2\mbox{\boldmath $t$}^2}{a^2b^2},
% $$
% tj.
% $$
% \left(\mbox{\boldmath $t$}^2-1\right) \left(\mbox{\boldmath $t$}^2w^2-a^2\right) = 0.
% $$
In this section we show that it is possible to unify these two cases$:$
\begin{equation}
\label{Fq}
\mathcal{F}{\mbox{\scriptsize \boldmath $q$}}: \quad 2\mbox{\boldmath $q$}^2wxy^2+\mbox{\boldmath $q$}^2b^2x^2
+(a^2+\mbox{\boldmath $q$}^4w^2)y^2-a^2b^2\mbox{\boldmath $q$}^2=0,
\end{equation}
using the parameter $\mbox{\boldmath $q$}$ given by the formula (\ref{QUu}). For the curve (\ref{Fq}), the following$:$
\begin{equation}
\label{Eq}
\mathcal{F}_{\mbox{\scriptsize \boldmath $q$}}
\equiv
\mathcal{F}
\end{equation}
holds if and only if
\begin{equation}
\frac{2 \mbox{\boldmath $q$}^2 w}{2 w}
=
\frac{\mbox{\boldmath $q$}^2 b^2}{b^2}
=
\frac{a^2 + \mbox{\boldmath $q$}^4w^2}{a^2+w^2}
=
\frac{a^2b^2\mbox{\boldmath $q$}^2}{a^2b^2}
\quad \Longleftrightarrow \quad
\left(\mbox{\boldmath $q$}^2-1\right) \left(\mbox{\boldmath $q$}^2w^2-a^2\right) = 0,
\end{equation}
i.e. the parameter $\mbox{\boldmath $q$}$ holds as defined in (\ref{QUu}).
Furthermore,
\begin{equation}
\label{UnionFq}
\mathcal{F}_{\mbox{\scriptsize \boldmath $q$}}
=
\mathcal{F}_{\mbox{\scriptsize \boldmath $q$},\, egg}
\cup
\mathcal{F}_{\mbox{\scriptsize \boldmath $q$},\, hyp},
\end{equation}
where $\mathcal{F}_{\mbox{\scriptsize \boldmath $q$},\, egg}$ is the egg-shaped part of $\mathcal{F}_{\mbox{\scriptsize \boldmath $q$}}$ over $[-a,a]$
and $\mathcal{F}_{\mbox{\scriptsize \boldmath $q$},\, hyp}$ is the hyperbolic-shaped part of $\mathcal{F}_{\mbox{\scriptsize \boldmath $q$}}$ over $(-\infty,\gamma{\mbox{\scriptsize \boldmath $q$}})$,
for $\gamma_{\mbox{\scriptsize \boldmath $q$}} = \gamma $.

\newpage\noindent
From the geometric point of view, H\" ugelsch\" affer's construction of the egg-shaped part of the curve is defined
using two non-concentric circles as considered in \cite{FH 1944}, \cite{Ferreol 2009}, \cite{MP 2010}, \cite{moNG 2010a}, \cite{Schmidbauer 1948} and \cite{Schmidbauer 1949}.
An analogous construction for $\mathcal{F}_{\mbox{\scriptsize \boldmath $q$},\, egg}$ is defined using the circles
$\mathcal{K}_{1}$ and $\mathcal{K}_{2}$ (see Fig. \ref{Fig.2}) given by the parametric equations
\begin{equation}
% \label{K1}
\mathcal{K}_{1}: \quad
\left\{
\begin{array}{l}
x_1(t) = a \cos t,  \\[1.0 ex]
y_1(t) = a \sin t
\end{array}
\right.
% \end{equation}
\qquad
\mbox{and}
\qquad
% \begin{equation}
% \label{K2}
\mathcal{K}_{2}: \quad
\left\{
\begin{array}{l}
x_2(t) = -\mbox{\boldmath $q$}^2 w + \mbox{\boldmath $q$} \, b \cos t,  \\[1.0 ex]
y_2(t) = \mbox{\boldmath $q$} \, b \sin t
\end{array}
\right.
\end{equation}
for $t \in [0, 2 \pi]$, such that the points
$
P_t = {\big (}x(t), y(t){\big )} \in \mathcal{F}_{\mbox{\scriptsize \boldmath $q$},\, egg}
$
have the following coordinates
\begin{equation}
\label{11}
\mathcal{F}_{\mbox{\scriptsize \boldmath $q$},\, egg}: \quad
\left\{
\begin{array}{l}
x(t) = -\mbox{\boldmath $q$}^2 w \sin^2 t + \cos t \sqrt{a^2 - \mbox{\boldmath $q$}^4 w^2 \sin^2 t},  \\[1.0 ex]
y(t) = \mbox{\boldmath $q$} \, b \sin t,
\end{array}
\right.
\end{equation}
for $t \in [0, 2 \pi]$.

%% FIGURE 2 %%%%%%%%%%%%%%%%%%%%%%%%%%%%%%%%%%%%%%%%%%%%%%%%%%%%%%%%%%%%%%%%%%%%%%%%%%%%%%%%%%%%%%%%

The parametrization (\ref{11}), for $0 \leq t \leq \pi$, determines the upper part
$\mathcal{F}_{\mbox{\scriptsize \boldmath $q$},\, egg}^{+}$
of curve $\mathcal{F}_{\mbox{\scriptsize \boldmath $q$},\, egg}$ from point $P_0 = (a,0)$
to point $P_{\pi} = (-a,0)$; while, for $\pi \leq t \leq 2 \pi$, the lower part
$\mathcal{F}_{\mbox{\scriptsize \boldmath $q$},\, egg}^{-}$
of curve $\mathcal{F}_{\mbox{\scriptsize \boldmath $q$},\, egg}$ is determined from point $P_{\pi} = (-a,0)$
to point $P_{2\pi} = P_{0} = (a,0)$. Let us note that $(\ref{UuN})$ holds for the abscissa
of point $P_{\pi/2}= (-\mbox{\boldmath $q$}^2 w,\mbox{\boldmath $q$} \, b)$.

% Za takav izbor parametarizacije $x(t)$ i $y(t)$ su neprekidno diferencijabilne funkcije nad $(0,\pi)$.

\bigskip
\noindent
The upper part $\mathcal{F}_{\mbox{\scriptsize \boldmath $q$},\, egg}^{+}$ can be further decomposed into an union of two disjunct portions
\begin{equation}
\mathcal{F}_{\mbox{\scriptsize \boldmath $q$},\, egg}^{+}
=
\ell_1 \cup \ell_2,
\end{equation}
such that for $0 \!\leq\! t \!<\! \pi/2$ we get $\ell_2$ and for $\pi/2 \!\leq\! t \!\leq\! \pi$ we get $\ell_1$ (see Fig. \ref{Fig.2}).

Let us note that in the paper \cite{AADM 2023}, surfaces $\mathcal{A}_1$ and $\mathcal{A}_2$ were considered
for the egg-shaped part $\mathcal{F}_{egg}$ of the H\" ugelsch\" affer curve from (\ref{UnionF})
(i.e. $\mathcal{F}_{\mbox{\scriptsize \boldmath $q$},\, egg}$ from (\ref{UnionFq}) because (\ref{Eq}) holds when $\mbox{\boldmath $q$} = 1$)
wherein the arc $\ell_1$ is part of the boundary of surface $\mathcal{A}_1$ and arc $\ell_2$ is part of the boundary of surface $\mathcal{A}_2$.

\section{The area of curve $\mathcal{F}_{\mbox{\scriptsize \boldmath $q$},\, egg}$}

\bigskip
In this paper, we give a new formula for calculating the area of the surface bound by the curve
$\mathcal{F}_{\mbox{\scriptsize \boldmath $q$},\, egg}$ using the elliptic integral of the first kind, \cite{Ryzhik 2015}:
\begin{equation*}
\mbox{\boldmath $K$}(k) \,=\! \int\limits_{0}^{\pi/2}  \!\frac{d \theta}{\sqrt{1-k^2\sin^2\theta}}, \;\;\;\; 0 \leq k^2 < 1,
\end{equation*}
and the elliptic integral of the second kind,  \cite{Ryzhik 2015}:
\begin{equation*}
\mbox{\boldmath $E$}(k) \,=\!\! \int\limits_{0}^{\pi/2} \!\sqrt{1 - k^2 \sin^2 \theta } \; d \theta, \;\;\;\; 0 \leq k^2 \leq 1.
\end{equation*}
Let it be noted that for $\mbox{\boldmath $K$}$ and $\mbox{\boldmath $E$}$, it holds that
\begin{equation*}
\mbox{\boldmath $K$}(k) = F\left(\frac{\pi}{2},k\right)
\qquad\mbox{and}\qquad
\mbox{\boldmath $E$}(k) = E\left(\frac{\pi}{2},k\right)
\end{equation*}
where $F(\theta,k)$ and $E(\theta,k)$ are elliptic integrals defined in \cite{Ryzhik 2015} with the formulae 8.112/1 and 8.112/2.
Let us also note that inequalities for these elliptic integrals were given recently in \cite{W-DJ 2023}.

\bigskip
\noindent
If the values of the parameter $t$ from $\pi$ to $0$ are considered, then the points $(x(t), y(t))$
are located on $\mathcal{F}^{+}_{\mbox{\scriptsize \boldmath $q$},\, egg}$, from $P_{\pi} = (-a,0)$ to $P_0 = (a,0)$.
In view of this, let $\mathcal{A}_{\mbox{\scriptsize \boldmath $q$},\,egg}$ denote the area of the surface bounded by the curve $\mathcal{F}_{\mbox{\scriptsize \boldmath $q$},\, egg}$ (see (\ref{11})), as given by the integral
\begin{equation}
\mathcal{A}_{\mbox{\scriptsize \boldmath $q$},\,egg}
=
2\!\!
\displaystyle \int\limits_{\pi}^{0}{ \!\! y(t) \, x'(t) \, dt}
=
-
2\!\!
\displaystyle \int\limits_{0}^{\pi}{ \!\! y(t) \, x'(t) \, dt}\,.
\end{equation}
Let
\begin{equation}
\mathcal{A}_{\mbox{\scriptsize \boldmath $q$},\,egg}
=
\mathcal{A}_{\mbox{\scriptsize \boldmath $q$},\,2}
+
\mathcal{A}_{\mbox{\scriptsize \boldmath $q$},\,1},
\end{equation}
where
\begin{equation}
\label{povrs_2}
\mathcal{A}_{\mbox{\scriptsize \boldmath$q$},\,2}
=
-2\!\!
\displaystyle \int\limits_{0}^{\pi/2}{ \!\! y(t) \, x'(t) \, dt}
\end{equation}
and
\begin{equation}
\label{povrs_1}
\mathcal{A}_{\mbox{\scriptsize \boldmath $q$},\,1}
=
-2\!\!
\displaystyle \int\limits_{\pi/2}^{\pi}{ \!\! y(t) \, x'(t) \, dt}\,.
\end{equation}

\medskip
\noindent
{\bf 1.} Firstly, $\mathcal{A}_{\mbox{\scriptsize \boldmath$q$},\,2}$ is calculated.
Substituting $x(t)$ and $y(t)$ from (\ref{11}) in (\ref{povrs_2}), we get
\begin{equation*}
\begin{array}{rcl}
\mathcal{A}_{\mbox{\scriptsize \boldmath $q$},\,2}
\!\!&\!\!=\!\!&\!\!
2\!\!
\displaystyle \int\limits_{0}^{\pi/2}{ \!\! \mbox{\boldmath $q$} \, b \sin t \,
{\Bigg (}
\, 2 \, \mbox{\boldmath $q$}^2 w \sin t \cos t
+
\sin t \sqrt{a^2 - \mbox{\boldmath $q$}^4 w^2 \sin^2 t}
+
\mbox{\small $\dfrac{\mbox{\boldmath $q$}^4 w^2 \sin t \cos^2 t}{\sqrt{a^2 - \mbox{\boldmath $q$}^4 w^2 \sin^2 t}}$}
\,{\Bigg )} \, dt}\,.
\end{array}
\end{equation*}
Let
\begin{equation}
I_1
=
\displaystyle \int\limits_{0}^{\pi/2}{ \!\! \sin^2 t \cos t \, dt},
\end{equation}
\begin{equation}
I_2
=
\displaystyle \int\limits_{0}^{\pi/2}{ \!\! \sin^2 t \sqrt{1 - \mbox{\small $\dfrac{\mbox{\boldmath $q$}^4 w^2}{a^2}$} \sin^2 t} \, dt},
\end{equation}
\begin{equation}
I_3
=
\displaystyle \int\limits_{0}^{\pi/2}{ \!\! \mbox{\small $\dfrac{\sin^2 t \cos^2 t}{\sqrt{1 - \mbox{\small $\dfrac{\mbox{\boldmath $q$}^4 w^2}{a^2}$} \sin^2 t}}$} \, dt}.
\end{equation}
Then
\begin{equation}
\mathcal{A}_{\mbox{\scriptsize \boldmath $q$},\,2}
=
4 \, \mbox{\boldmath $q$}^3 w \, b \, I_1 + 2 \, \mbox{\boldmath $q$} \, a \, b \, I_2 + 2 \, \frac{\mbox{\boldmath $q$}^5 b \, w^2}{a} \, I_3 \,.
\end{equation}
It holds that
\begin{equation}
I_1
=
\dfrac{1}{3}
\end{equation}
For $I_2$, formula 2.583/4 from \cite{Ryzhik 2015} is used and it follows that
\begin{equation}
\begin{array}{rcl}
I_2
\!\!&\!\!=\!\!&\!\!
\left(-
\dfrac{\sin t \cos t}{3}\sqrt{1 - k^2 \sin^2 t}
+
\dfrac{1-k^2}{3k^2} F(t,k)
+
\dfrac{2 k^2 - 1}{3 k^2} E(t,k)
\right) \!\!\mathop{\mbox{\LARGE $|$}}\limits_{0}^{\pi/2}                           \\[2.0 ex]
\!\!&\!\!=\!\!&\!\!
\dfrac{1-k^2}{3k^2} F(\frac{\pi}{2},k)
+
\dfrac{2 k^2 - 1}{3 k^2} E(\frac{\pi}{2},k)
-
\dfrac{1-k^2}{3k^2} F(0,k)
-
\dfrac{2 k^2 - 1}{3 k^2} E(0,k)                           \\[2.0 ex]
\!\!&\!\!=\!\!&\!\!
\dfrac{1-k^2}{3k^2} \mbox{\boldmath $K$}(k)
+
\dfrac{2 k^2 - 1}{3 k^2} \mbox{\boldmath $E$}(k),
\end{array}
\end{equation}
where $k^2 = \dfrac{\mbox{\boldmath $q$}^4 w^2}{a^2}$ and it is true that $0 < k^2 < 1$.

\medskip
\noindent
For $I_3$ formula 2.584/13 from \cite{Ryzhik 2015} is used and it follows that
\begin{equation}
\begin{array}{rcl}
I_3
\!\!&\!\!=\!\!&\!\!
\left(
-
\dfrac{\sin t \cos t}{3k^2}\sqrt{1 - k^2 \sin^2 t}
+
\dfrac{2k^2-2}{3k^4} F(t,k)
+
\dfrac{2 - k^2}{3 k^4} E(t,k)
\right)\!\!\mathop{\mbox{\LARGE $|$}}\limits_{0}^{\pi/2}                           \\[2.0 ex]
\!\!&\!\!=\!\!&\!\!
\dfrac{2k^2-2}{3k^4} F(\frac{\pi}{2},k)
+
\dfrac{2 - k^2}{3 k^4} E(\frac{\pi}{2},k)
-
\dfrac{2k^2-2}{3k^4} F(0,k)
-
\dfrac{2 - k^2}{3 k^4} E(0,k)                           \\[2.0 ex]
\!\!&\!\!=\!\!&\!\!
\dfrac{2k^2-2}{3k^4}  \mbox{\boldmath $K$}(k)
+
\dfrac{2 - k^2}{3 k^4} \mbox{\boldmath $E$}(k),
\end{array}
\end{equation}
where $k^2 = \dfrac{\mbox{\boldmath $q$}^4 w^2}{a^2}$ and it is true that $0 < k^2 < 1$.
Based on the previous, it holds that
\begin{equation}
\begin{array}{rcl}
\mathcal{A}_{\mbox{\scriptsize \boldmath $q$},\,2}
\!\!&\!\!=\!\!&\!\!
4 \, \mbox{\boldmath $q$}^3 w \, b \, I_1
+ 2 \, \mbox{\boldmath $q$} \, a \, b \, I_2
+ 2 \, \dfrac{\mbox{\boldmath $q$}^5 b \, w^2}{a} \, I_3                           \\[2.0 ex]
\!\!&\!\!=\!\!&\!\!
\dfrac{4}{3} \, \mbox{\boldmath $q$}^3 w \, b
+ 2 \, \mbox{\boldmath $q$} \, a \, b \,
\left(
\dfrac{1-k^2}{3k^2} \mbox{\boldmath $K$}(k)
+
\dfrac{2 k^2 - 1}{3 k^2} \mbox{\boldmath $E$}(k)
\right)+                                                                      \\[2.0 ex]
\!\!&\!\!\!\!&\!\!+\;
2 \, \dfrac{\mbox{\boldmath $q$}^5 b \, w^2}{a} \,
\left(
\dfrac{2k^2-2}{3k^4}  \mbox{\boldmath $K$}(k)
+
\dfrac{2 - k^2}{3 k^4} \mbox{\boldmath $E$}(k)
\right)                                                                         \\[2.0 ex]
\!\!&\!\!=\!\!&\!\!
\dfrac{2}{3}\, a \, b \, \mbox{\boldmath $q$}
\left( \left(1-\dfrac{1}{k^2}\right)\mbox{\boldmath $K$}(k)
+
\left(1+\dfrac{1}{k^2}\right) \mbox{\boldmath $E$}(k)
+
2 k
\right).
\end{array}
\end{equation}

\bigskip
\noindent
{\bf 2.} To calculate $\mathcal{A}_{\mbox{\scriptsize \boldmath $q$},\,1}$,
substituting $x(t)$ and $y(t)$ from (\ref{11}) in (\ref{povrs_1}), we get
\begin{equation*}
\begin{array}{rcl}
\mathcal{A}_{\mbox{\scriptsize \boldmath $q$},\,1}
\!\!&\!\!=\!\!&\!\!
2\!\!
\displaystyle \int\limits_{\pi/2}^{\pi}{ \!\! \mbox{\boldmath $q$} \, b \sin t \,
{\Bigg (}
\, 2 \, \mbox{\boldmath $q$}^2 w \sin t \cos t
+ \sin t \sqrt{a^2 - \mbox{\boldmath $q$}^4 w^2 \sin^2 t}
+ \mbox{\small $\dfrac{\mbox{\boldmath $q$}^4 w^2 \sin t \cos^2 t}{\sqrt{a^2 - \mbox{\boldmath $q$}^4 w^2 \sin^2 t}}$}
\,{\Bigg )} \, dt}\,.
\end{array}
\end{equation*}
Let
\begin{equation}
J_1
=
\displaystyle \int\limits_{\pi/2}^{\pi}{ \!\! \sin^2 t \cos t \, dt},
\end{equation}
\begin{equation}
J_2
=
\displaystyle \int\limits_{\pi/2}^{\pi}{ \!\! \sin^2 t \sqrt{1 - \mbox{\small $\dfrac{\mbox{\boldmath $q$}^4 w^2}{a^2}$} \sin^2 t} \, dt},
\end{equation}
\begin{equation}
J_3
=
\displaystyle \int\limits_{\pi/2}^{\pi}{ \!\! \mbox{\small $\dfrac{\sin^2 t \cos^2 t}{\sqrt{1 - \mbox{\small $\dfrac{\mbox{\boldmath $q$}^4 w^2}{a^2}$} \sin^2 t}}$} \, dt}.
\end{equation}
Then
\begin{equation}
\mathcal{A}_{\mbox{\scriptsize \boldmath $q$},\,1}
=
4 \, \mbox{\boldmath $q$}^3 w \, b \, J_1
+ 2 \, \mbox{\boldmath $q$} \, a \, b \, J_2
+ 2 \, \frac{\mbox{\boldmath $q$}^5 b \, w^2}{a} \, J_3 \,.
\end{equation}
It holds that
\begin{equation}
J_1
=
- \dfrac{1}{3},\;\;
J_2
=
\displaystyle \int\limits_{0}^{\pi/2}{ \!\! \cos^2 t
\sqrt{1 - \mbox{\small $\dfrac{\mbox{\boldmath $q$}^4 w^2}{a^2}$} \cos^2 t} \, dt},
\;\;
J_3
=
\displaystyle \int\limits_{0}^{\pi/2}{ \!\! \mbox{\small $\dfrac{\cos^2 t \sin^2 t}{
\sqrt{1 - \mbox{\small $\dfrac{\mbox{\boldmath $q$}^4 w^2}{a^2}$} \cos^2 t}}$} \, dt}.
\end{equation}
Analogously to the previous
\begin{equation}
J_2
\!=\!
\dfrac{1\!-\!k^2}{3 k^2} \mbox{\boldmath $K$}(k)
+
\dfrac{2 k^2\!-\!1}{3 k^2} \mbox{\boldmath $E$}(k)
\!=\!
I_2,
\;\;
J_3
\!=\!
\dfrac{2 k^2\!-\!2}{3 k^4} \mbox{\boldmath $K$}(k)
+
\dfrac{2\!-\!k^2}{3 k^4} \mbox{\boldmath $E$}(k)
\!=\!
I_3,
\end{equation}
where $k^2 = \dfrac{\mbox{\boldmath $q$}^4 w^2}{a^2}$ and it is true that $0 < k^2 < 1$.
Thus, we obtain that
\begin{equation}
\begin{array}{rcl}
\mathcal{A}_{\mbox{\scriptsize \boldmath $q$},\,1}
\!\!&\!\!=\!\!&\!\!
4 \, \mbox{\boldmath $q$}^3 w \, b \, J_1
+ 2 \, \mbox{\boldmath $q$} \, a \, b \, J_2
+ 2 \, \displaystyle\frac{\mbox{\small $q$}^5 b \, w^2}{a} \, J_3 \\[2.0 ex]
\!\!&\!\!=\!\!&\!\!
\dfrac{2}{3}\, a \, b \,\mbox{\boldmath $q$}
\left(
\left(1-\dfrac{1}{k^2}\right)\mbox{\boldmath $K$}(k)
+
\left(1+\dfrac{1}{k^2}\right)\mbox{\boldmath $E$}(k)
-
2 k
\right).
\end{array}
\end{equation}

\medskip
\noindent
In all, the following theorem has been proven.

\smallskip
\noindent
\begin{theorem}
For the area $\mathcal{A}_{\mbox{\scriptsize \boldmath$q$},\,egg}$
of curve $\mathcal{F}_{\mbox{\scriptsize \boldmath $q$},\, egg}$ it holds that$:$
\begin{equation}
\label{A1&A2}
\mathcal{A}_{\mbox{\scriptsize \boldmath $q$},\,egg}
=
\mathcal{A}_{\mbox{\scriptsize \boldmath $q$},\,2}+\mathcal{A}_{\mbox{\scriptsize \boldmath $q$},\,1}
=
\dfrac{4}{3}\, a \, b \,\mbox{\boldmath $q$}
\left(
\left(1-\frac{1}{k^2}\right)\mbox{\boldmath $K$}(k)
+
\left(1+\frac{1}{k^2}\right)\mbox{\boldmath $E$}(k)
\right)
\end{equation}
and
\begin{equation}
\label{prethodna}
\mathcal{A}_{\mbox{\scriptsize \boldmath $q$},\,2}-\mathcal{A}_{\mbox{\scriptsize \boldmath $q$},\,1}
=
\dfrac{8}{3} \, a \, b \,\mbox{\boldmath $q$} \, k
=
\dfrac{8}{3}\, w \, b \,\mbox{\boldmath $q$}^3
= % \mathop{=}\limits_{\,(\ref{QUu})\,}
\left\{
\begin{array}{ccc}
8\, w \, b /3    &:& w < a, \\[1.0 ex]
\dfrac{8}{3\, w ^2}\, b \, a^3   &:& w > a;
\end{array}
\right.
\end{equation}
where $a$, $b$, $w$ are H\" ugelsch\" affer curve parameters and $k = \dfrac{\mbox{\boldmath $q$}^2 w}{a}$.
\end{theorem}

\smallskip
\noindent
\begin{remark}
Let us note that, for cases $(I)$ and $(II)$, the equality (\ref{prethodna}) is derived in \cite{AADM 2023} (p. 185).
\end{remark}

Alongside the initial applications of H\" ugelsch\" affer curves in aero-engineering (see \cite{Boermans 2004}, \cite{FH 1944}), recently, there has been research on the
applications of these curves in: architecture and civil engineering (see \cite{MP 2010}, \cite{ICEGD 2011}); poultry industry, ornithology and bioengineering (see \cite{HK 2023}, \cite{EggCurve 2019}, \cite{EggCurve 2020}, \cite{Narushin 2020}, \cite{Narushin 2021}, \cite{Narushin 2022}, \cite{Narushin 2023}); traffic engineering (see \cite{PS 2023}) and
hydro-engineering (see \cite{WaterS 2021}, \cite{AADM 2023}, \cite{WaterER 2022}). To aid in the application of H\" ugelsch\" affer curves and the practical usage of the area formulae for these curves, we have developed the applet \cite{Applet 2022}.

\section{Taylor approximations of elliptic integrals}

\bigskip
\noindent
Let $f : (\alpha, \beta) \longrightarrow \RR $ be a real function. We state some definitions and characteristics according to \cite{MRL 2019}.
\begin{definition}
Let $T_n^{f,\,\alpha_+}(x)$ be a Taylor polynomial for function $f(x)$, of degree $n \in \NN_0$, in the right neighborhood of point $\alpha$.
If, for the real function $f : (\alpha, \beta) \longrightarrow \RR $, there exist finite limits
$f^{(i)}(\alpha_+) = \lim_{x\rightarrow \alpha_+}{f^{(i)}(x)} $, for $i \in\{0,1,...,n\}$ then
\begin{equation}
T_n^{f,\,\alpha_+}(x) = \sum_{i = 0}^{n} \frac{f^{(i)}(\alpha_+)}{i!}(x-\alpha)^{i}
\end{equation}
is the first Taylor approximation of function $f$ in the right neighborhood of point $\alpha$, for $n \in \NN_0$, where
\begin{equation}
R_n^{f,\,\alpha_+}(x) = f(x) - T_{n-1}^{f,\,\alpha_+}(x)
\end{equation}
is the remainder of the first Taylor approximation in the right neighborhood of point $\alpha$.
\end{definition}
\begin{definition}
Let $T_n^{f,\,\beta_-}(x)$ be a Taylor polynomial for function $f(x)$, of degree $n \in \NN_0$, in the left neighborhood of point $\beta$.
If, for the real function $f : (\alpha,\beta) \longrightarrow \RR $, there exist finite limits
$f^{(i)}(\beta_-) = \lim_{x\rightarrow \beta_-}{f^{(i)}(x)}$, for $i \in\{0,1,...,n\}$ then
\begin{equation}
T_n^{f,\,\beta_-}(x) = \sum_{i = 0}^{n} \frac{f^{(i)}(\beta_-)}{i!}(x-\beta)^{i}
\end{equation}
is the first Taylor approximation of function $f$ in the left neighborhood of point $\beta$, for $n \in \NN_0$, where
\begin{equation}
R_n^{f,\,\beta_-}(x) = f(x) - T_{n-1}^{f,\,\beta_-}(x)
\end{equation}
is the remainder of the first Taylor approximation in the left neighborhood of point $\beta$.
\end{definition}

\begin{definition}
For the polynomial of the form
\begin{equation}
\mathbbmsl{T}_n^{f;\,\alpha_+,\,\beta_-}(x) =
    \left\{
    \begin{array}{lll}
    T_{n-1}^{f,\,\alpha_+}(x) + \dfrac{(x-\alpha)^{n}}{(\beta-\alpha)^{n}}R_{n}^{f,\,\alpha_+}(\beta_-) & : & n \geq 1 \\[1.0 ex]
    f(\beta_-) & : & n = 0
    \end{array}
    \right.
\end{equation}
we say that it is the second Taylor approximation of function $f$ in the right neighborhood of point $\alpha$, for $n \in \NN_0$, while the polynomial
\begin{equation}
\mathbbmsl{T}_n^{f;\,\beta_-,\,\alpha_+}(x) =
    \left\{
    \begin{array}{lll}
    T_{n-1}^{f,\,\beta_-}(x) + \dfrac{(x-\beta)^{n}}{(\alpha-\beta)^{n}}R_{n}^{f,\,\beta_-}(\alpha_+) & : & n \geq 1 \\[1.0 ex]
    f(\alpha_+) & : & n = 0
    \end{array}
    \right.
\end{equation}
is the second Taylor approximation of function $f$ in the left neighborhood of point $\beta$, for $n \in \NN_0$.
\end{definition}
\begin{theorem}[Theorem WD] \label{Theorem WD}
Suppose that function $f$ is real over $(\alpha,\beta)$, i.e. $f : (\alpha,\beta) \longrightarrow \RR $
and let $n$ be a whole natural number such that $f^{(i)}(\alpha_+)$ and $f^{(i)}(\beta_-)$ exist, for $i \in \{0,1,...,n\}$. \\ [+1.75 ex]
Also suppose that $(-1)^{(n)} f^{(n)}(x)$ is increasing over $(\alpha,\beta)$.
Then, for each $x \in (\alpha,\beta)$ the following inequality holds:
\begin{equation}
\mathbbmsl{T}_n^{f;\,\beta_-,\,\alpha_+}(x) < f(x) < T_n^{f;\,\beta_-}(x)
\end{equation}
and supposing that $f^{(n)}(x)$ is increasing over $(\alpha,\beta)$, then for each $x \in (\alpha,\beta)$ it holds that:
\begin{equation}
\mathbbmsl{T}_n^{f;\,\alpha_+,\,\beta_-}(x) > f(x) > T_n^{f;\,\alpha+}(x)
\end{equation}
When the function $(-1)^{(n)} f^{(n)}(x)$ is decreasing over $(\alpha,\beta)$, or
when $f^{(n)}(x)$ is decreasing over $(\alpha,\beta)$ then for each $x \in (\alpha,\beta)$
the reverse inequalities hold.
\end{theorem}
The preceding theorem was proven by S. Wu and L. Debnath u \cite{WD 2009}.
In \cite{MRL 2019} and \cite{MRL 2022} some applications of this statement were considered within the Theory of analytical inequalities, see also papers
\cite{BagulMalesevic 2023} and \cite{ChaoMalesevic 2023}.

\bigskip
\noindent
{\bf \boldmath 3.1 Taylor approximations of $K$ and $E$ elliptic integrals} %
% %%%%%%%%%%%%%%%%%%%%%%%%%%%%%%%%%%%%%%%%%%%%%%%%%%%%%%%%%%%%%%%%%%%%%%%%%%%%%%%%%%%%%%%%%%%%%%%%%%%%%%%%%%%%%%%%%%%%%%%%%%%%%%%%%%%%%%%%%%%%%%%%%%%%%

\bigskip
%\noindent
Let us apply the previous consideration to elliptic integrals $\mbox{\boldmath $K$}(x)$ over $(-1,1)$
and $\mbox{\boldmath $E$}(x)$ over $[-1,1]$. We will use the well-known series expansions:
\begin{equation}
\label{K_series}
\mbox{\boldmath $K$}(x) = \frac{\pi}{2} \sum_{i=0}^{\infty}{\left( \frac{(2i-1)!!}{(2i)!!} \right)^2 \! x^{2i}},
\end{equation}
for $x \!\in\! (-1,1)$, see formula 8.113/1 in \cite{Ryzhik 2015} and
\begin{equation}
\label{E_series}
\mbox{\boldmath $E$}(x) =
\frac{\pi}{2} - \frac{\pi}{2} \sum_{i=1}^{\infty}{ \left( \frac{(2i-1)!!}{(2i)!!} \right)^2 \!\! \frac{1}{2i-1}} \, x^{2i},
\end{equation}
for $x \!\in\! [-1,1]$,  see formula 8.114/1 in \cite{Ryzhik 2015}.
Let us note that the series expansion (\ref{E_series}) has a radius of $1$ and that $x = 1$ can be included in the convergence domain,
and that the series expansion (\ref{K_series}) has a radius of $1$ and that $x = 1$ can not be included in the convergence domain.

\medskip
\noindent
For the elliptic integral $\mbox{\boldmath $K$}(x)$, based on the series expansion (\ref{K_series}), the first and second Taylor approximations
are obtained. The first Taylor approximation is:
\begin{equation}
T_n^{\mbox{\scriptsize \boldmath $K$},\,0}(x)
=
\frac{\pi}{2} \sum_{i=0}^{[n/2]}{\left( \frac{(2i-1)!!}{(2i)!!} \right)^2} \! x^{2i},
\end{equation}
where $x \!\in\! (-1,1)$. For a fixed $\beta \!\in\! (0,1)$, the second Taylor approximation in the right neighborhood of point
$\alpha = 0$ is:
\begin{equation}
\mathbbmsl{T}_n^{\mbox{\scriptsize \boldmath $K$};\,0,\beta}(x) =
    \left\{
    \begin{array}{lll}
    T_{n-1}^{\mbox{\scriptsize \boldmath $K$},\,0}(x) + \dfrac{x^n}{\beta^n}\left(\mbox{\boldmath $K$}(\beta) - T_{n-1}^{\mbox{\scriptsize \boldmath $K$},\,0}(\beta) \right) & : & n \geq 1 \\[1.0 ex]
    \mbox{\boldmath $K$}(\beta) & : & n = 0
    \end{array}
    \right.,
\end{equation}
where $x \!\in\! [-\beta,\beta]$.

\medskip
\noindent
For the elliptic integral $\mbox{\boldmath $E$}(x)$, based on the series expansion (\ref{E_series}), the first and second Taylor approximations
are obtained. The first Taylor approximation is:
\begin{equation}
T_n^{\mbox{\scriptsize \boldmath $E$},\,0}(x)
=
\frac{\pi}{2} - \frac{\pi}{2} \sum_{i=1}^{[n/2]}{\left( \frac{(2i-1)!!}{(2i)!!} \right)^2 \!\! \frac{1}{2i-1}} \, x^{2i},
\end{equation}
where $x \!\in\! [-1,1]$. For a fixed $\beta \!\in\! (0,1]$, the second Taylor approximation in the right neighborhood of point
$\alpha = 0$ is:
\begin{equation}
\mathbbmsl{T}_n^{\mbox{\scriptsize \boldmath $E$};\,0,\,\beta}(x) =
    \left\{
    \begin{array}{lll}
    T_{n-1}^{\mbox{\scriptsize \boldmath $E$},\,0}(x) + \dfrac{x^n}{\beta^n}\left(\mbox{\boldmath $E$}(\beta) - T_{n-1}^{\mbox{\scriptsize \boldmath $E$},\,0}(\beta) \right) & : & n \geq 1 \\[1.0 ex]
    \mbox{\boldmath $E$}(\beta) & : & n = 0
    \end{array}
    \right.,
\end{equation}
where $x \!\in\! [-\beta,\beta]$.

\medskip
\noindent
According to \cite{MRL 2019}, the following inequalities hold for the elliptic integral $\mbox{\boldmath $K$}(x):$
\begin{equation}
\label{Nejednakosti_za_K}
\begin{array}{c}
\dfrac{\pi}{2}
\!=\!
T_{0}^{\mbox{\scriptsize \boldmath $K$},\,0}(x) \!=\! T_{1}^{\mbox{\scriptsize \boldmath $K$},\,0}(x) \leq
T_{2}^{\mbox{\scriptsize \boldmath $K$},\,0}(x) \!=\! T_{3}^{\mbox{\scriptsize \boldmath $K$},\,0}(x) \leq
\ldots \leq
T_{2i}^{\mbox{\scriptsize \boldmath $K$},\,0}(x) \!=\! T_{2i+1}^{\mbox{\scriptsize \boldmath $K$},\,0}(x) \leq
\ldots                                                                          \\[1.5 ex]
\leq \mbox{\boldmath $K$}(x) \leq                                                                        \\[1.5 ex]
\ldots \leq
\mathbbmsl{T}_j^{\mbox{\scriptsize \boldmath $K$};\,0,\,\beta}(x) \leq
\ldots               \leq
\mathbbmsl{T}_3^{\mbox{\scriptsize \boldmath $K$};\,0,\,\beta}(x) \leq
\mathbbmsl{T}_2^{\mbox{\scriptsize \boldmath $K$};\,0,\,\beta}(x) \leq
\mathbbmsl{T}_1^{\mbox{\scriptsize \boldmath $K$};\,0,\,\beta}(x) \leq
\mathbbmsl{T}_0^{\mbox{\scriptsize \boldmath $K$};\,0,\,\beta}(x)
\!=\!
\mbox{\boldmath $K$}(\beta)
\end{array}
\end{equation}
for a fixed $\beta \!\in\! (0,1)$ and an arbitrary $x \!\in\! [-\beta,\beta]$. Furthermore, according to \cite{MRL 2019}
the following inequalities hold for the elliptic integral $\mbox{\boldmath $E$}(x):$
\begin{equation}
\label{Nejednakosti_za_E}
\begin{array}{c}
\dfrac{\pi}{2}
\!=\!
T_{0}^{\mbox{\scriptsize \boldmath $E$},\,0}(x) \!=\! T_{1}^{\mbox{\scriptsize \boldmath $E$},\,0}(x) \geq
T_{2}^{\mbox{\scriptsize \boldmath $E$},\,0}(x) \!=\! T_{3}^{\mbox{\scriptsize \boldmath $E$},\,0}(x) \geq
\ldots \geq
T_{2i}^{\mbox{\scriptsize \boldmath $E$},\,0}(x) \!=\! T_{2i+1}^{\mbox{\scriptsize \boldmath $E$},\,0}(x) \geq
\ldots                                                                          \\[2 ex]
\geq \mbox{\boldmath $E$}(x) \geq                                                                        \\[2 ex]
\ldots \geq
\mathbbmsl{T}_j^{\mbox{\scriptsize \boldmath $E$};\,0,\,\beta}(x) \geq
\ldots               \geq
\mathbbmsl{T}_3^{\mbox{\scriptsize \boldmath $E$};\,0,\,\beta}(x) \geq
\mathbbmsl{T}_2^{\mbox{\scriptsize \boldmath $E$};\,0,\,\beta}(x) \geq
\mathbbmsl{T}_1^{\mbox{\scriptsize \boldmath $E$};\,0,\,\beta}(x) \geq
\mathbbmsl{T}_0^{\mbox{\scriptsize \boldmath $E$};\,0,\,\beta}(x)
\!=\!
\mbox{\boldmath $E$}(\beta)
\end{array}
\end{equation}
for a fixed $\beta \!\in\! (0,1]$ and an arbitrary $x \!\in\! [-\beta, \beta]$.

\noindent
In table 1, we list the explicit forms of the polynomials $T_j^{{\mbox{\scriptsize \boldmath $K$}}^{\vphantom{2^2}},\,0}(x)$
and $\mbox{$\mathbbmsl{T}$}_j^{\mbox{\scriptsize \boldmath $K$};\,0,\,\beta}(x)$,
for $\beta~\!\in\!~(0,1)$ and $j = 0, 1, \ldots, 10;$
while, in table 2, we list the explicit forms of the polynomials $T_j^{{\mbox{\scriptsize \boldmath $E$}}^{\vphantom{2^2}},\,0}(x)$
and $\mbox{$\mathbbmsl{T}$}_j^{\mbox{\scriptsize \boldmath $E$};\,0,\,1}(x)$,
for $j = 0, 1, \ldots, 10.$

%\begin{landscape}
%\medskip
\begin{equation*}
\begin{array}{|c|l|l|} \hline
    &
    &                                                                           \\[0.15 ex]
 j & T_j^{{\mbox{\scriptsize \boldmath $K$}}^{\vphantom{2^2}},\,0}(x)                                            & \mbox{$\mathbbmsl{T}$}_j^{\mbox{\scriptsize \boldmath $K$};\,0\,,\beta}(x) \;\mbox{and}\; C_{\mbox{\tiny \boldmath $K$}} = \mbox{\boldmath $K$}(\beta) \;\mbox{for}\; \beta \!\in\! (0,1)                                                       \\[1.5 ex]  \hline \hline
     &
    &                                                                           \\[0.15 ex]
 0 & \nnfrac{{\pi}^{\vphantom{2}}}{2}                                             & C_{\mbox{\tiny \boldmath $K$}}                                                                                 \\[1.5 ex]
 1 & \nnfrac{\pi}{2}                                                              & T_0^{\mbox{\scriptsize \boldmath $K$};\,0}(x)\!+\!\left(\nnfrac{C_{\mbox{\tiny \boldmath $K$}}}{\beta}\!-\!\nnfrac{\pi}{2\beta}\right)x                              \\[1.5 ex]
 2 & \nnfrac{\pi}{2}\!+\!\nnfrac{\pi}{8}x^2                                           & T_1^{\mbox{\scriptsize \boldmath $K$};\,0}(x)\!+\!\left(\nnfrac{C_{\mbox{\tiny \boldmath $K$}}}{\beta^2}\!-\!\nnfrac{\pi}{2\beta^2}\right)x^2                        \\[1.5 ex]
 3 & \nnfrac{\pi}{2}\!+\!\nnfrac{\pi}{8}x^2                                           & T_2^{\mbox{\scriptsize \boldmath $K$};\,0}(x)\!+\!\left(\nnfrac{C_{\mbox{\tiny \boldmath $K$}}}{\beta^3}\!-\!\nnfrac{\pi}{2\beta^3}\!-\!\nnfrac{\pi}{8\beta}\right)x^3       \\[1.5 ex]
 4 & \nnfrac{\pi}{2}\!+\!\nnfrac{\pi}{8}x^2\!+\!\nnfrac{9\pi}{128}x^4                     & T_3^{\mbox{\scriptsize \boldmath $K$};\,0}(x)\!+\!\left(\nnfrac{C_{\mbox{\tiny \boldmath $K$}}}{\beta^4}\!-\!\nnfrac{\pi}{2\beta^4}\!-\!\nnfrac{\pi}{8\beta^2}\right)x^4     \\[1.5 ex]
 5 & \nnfrac{\pi}{2}\!+\!\nnfrac{\pi}{8}x^2\!+\!\nnfrac{9\pi}{128}x^4
   & T_4^{\mbox{\scriptsize \boldmath $K$};\,0}(x)\!+\!\left(\nnfrac{C_{\mbox{\tiny \boldmath $K$}}}{\beta^5}\!-\!\nnfrac{\pi}{2\beta^5}\!-\!\nnfrac{\pi}{8\beta^3}\!-\!\nnfrac{9\pi}{128\beta}\right)x^5   \\[1.5 ex]
 6 & \nnfrac{\pi}{2}\!+\!\nnfrac{\pi}{8}x^2\!+\!\nnfrac{9\pi}{128}x^4\!+\!\nnfrac{25\pi}{512}x^6
   & T_5^{\mbox{\scriptsize \boldmath $K$};\,0}(x)\!+\!\left(\nnfrac{C_{\mbox{\tiny \boldmath $K$}}}{\beta^6}\!-\!\nnfrac{\pi}{2\beta^6}\!-\!\nnfrac{\pi}{8\beta^4}\!-\!\nnfrac{9\pi}{128\beta^2}\right)x^6  \\[1.5 ex]
 7 & \nnfrac{\pi}{2}\!+\!\nnfrac{\pi}{8}x^2\!+\!\nnfrac{9\pi}{128}x^4\!+\!\nnfrac{25\pi}{512}x^6
   & T_6^{\mbox{\scriptsize \boldmath $K$};\,0}(x)\!+\!\left(\nnfrac{C_{\mbox{\tiny \boldmath $K$}}}{\beta^7}\!-\!\nnfrac{\pi}{2\beta^7}\!-\!\nnfrac{\pi}{8\beta^5}\!-\!\nnfrac{9\pi}{128\beta^3}\!-\!\nnfrac{25\pi}{512\beta}\right)x^7 \\[1.5 ex]
 8 & \nnfrac{\pi}{2}\!+\!\nnfrac{\pi}{8}x^2\!+\!\nnfrac{9\pi}{128}x^4\!+\!\nnfrac{25\pi}{512}x^6
   \!\!+\!\! \nnfrac{1225\pi}{32768}x^8
   & T_7^{\mbox{\scriptsize \boldmath $K$};\,0}(x)\!+\!\left(\nnfrac{C_{\mbox{\tiny \boldmath $K$}}}{\beta^8}\!-\!\nnfrac{\pi}{2\beta^8}\!-\!\nnfrac{\pi}{8\beta^6}\!-\!\nnfrac{9\pi}{128\beta^4}\!-\!\nnfrac{25\pi}{512\beta^2}\right)x^8 \\[1.5 ex]
 9 & \nnfrac{\pi}{2}\!+\!\nnfrac{\pi}{8}x^2\!+\!\nnfrac{9\pi}{128}x^4\!+\!\nnfrac{25\pi}{512}x^6
   \!\!+\!\! \nnfrac{1225\pi}{32768}x^8
   & T_8^{\mbox{\scriptsize \boldmath $K$};\,0}(x)\!+\!\left(\nnfrac{C_{\mbox{\tiny \boldmath $K$}}}{\beta^9}\!-\!\nnfrac{\pi}{2\beta^9}\!-\!\nnfrac{\pi}{8\beta^7}\!-\!\nnfrac{9\pi}{128\beta^5}\!-\!\nnfrac{25\pi}{512\beta^3}\!-\!\nnfrac{1225\pi}{32768\beta}\right)x^9  \\[1.5 ex]
10 & \nnfrac{\pi}{2}\!+\!\nnfrac{\pi}{8}x^2\!+\!\nnfrac{9\pi}{128}x^4\!+\!\nnfrac{25\pi}{512}x^6
   \!\!+\!\! \nnfrac{1225\pi}{32768}x^8
   \!\!+\!\! \nnfrac{3969\pi}{131072}x^{10} \,\,
   & T_9^{\mbox{\scriptsize \boldmath $K$};\,0}(x)\!+\!\left(\nnfrac{C_{\mbox{\tiny \boldmath $K$}}}{\beta^{10}}\!-\!\nnfrac{\pi}{2\beta^{10}}\!-\!\nnfrac{\pi}{8\beta^8}\!-\!\nnfrac{9\pi}{128\beta^6}\!-\!\nnfrac{25\pi}{512\beta^4}\!-\!\nnfrac{1225\pi}{32768\beta^2}\right)x^{10}  \\[1.5 ex]

\hline
\end{array}
\end{equation*}

\vspace*{0.00 mm}

\centerline{\hspace*{-15.0 mm} Table 1.}
%\end{landscape}

%\begin{landscape}
%\medskip
\begin{equation*}
\begin{array}{|c|l|l|} \hline
    &
    &                                                                           \\[0.05 ex]
 j & T_j^{{\mbox{\scriptsize \boldmath $E$}}^{\vphantom{2^2}},\,0}(x)                                            & \mbox{$\mathbbmsl{T}$}_j^{\mbox{\scriptsize \boldmath $E$};\,0,\,1}(x)\;\mbox{and}\;\beta = 1                                                               \\[1.5 ex]  \hline \hline
     &
    &                                                                           \\[0.05 ex]
 0 & \nnfrac{{\pi}^{\vphantom{2}}}{2}                                             & 1                                                                            \\[1.75 ex]
 1 & \nnfrac{\pi}{2}                                                              & T_0^{\mbox{\scriptsize \boldmath $E$},\,0}(x)+\left(1-\nnfrac{\pi}{2}\right)x                                         \\[1.75 ex]
 2 & \nnfrac{\pi}{2}-\nnfrac{\pi}{8}x^2                                           & T_1^{\mbox{\scriptsize \boldmath $E$},\,0}(x)+\left(1-\nnfrac{\pi}{2}\right)x^2                                       \\[1.75 ex]
 3 & \nnfrac{\pi}{2}-\nnfrac{\pi}{8}x^2                                           & T_2^{\mbox{\scriptsize \boldmath $E$},\,0}(x)+\left(1-\nnfrac{3\pi}{8}\right)x^3                                       \\[1.75 ex]
 4 & \nnfrac{\pi}{2}-\nnfrac{\pi}{8}x^2-\nnfrac{3\pi}{128}x^4                     & T_3^{\mbox{\scriptsize \boldmath $E$},\,0}(x)+\left(1-\nnfrac{3\pi}{8}\right)x^4                                       \\[1.75 ex]
 5 & \nnfrac{\pi}{2}-\nnfrac{\pi}{8}x^2-\nnfrac{3\pi}{128}x^4                             & T_4^{\mbox{\scriptsize \boldmath $E$},\,0}(x)+\left(1-\nnfrac{45\pi}{128}\right)x^5                            \\[1.75 ex]
 6 & \nnfrac{\pi}{2}-\nnfrac{\pi}{8}x^2-\nnfrac{3\pi}{128}x^4-\nnfrac{5\pi}{512}x^6       & T_5^{\mbox{\scriptsize \boldmath $E$},\,0}(x)+\left(1-\nnfrac{45\pi}{128}\right)x^6                            \\[1.75 ex]
 7 & \nnfrac{\pi}{2}-\nnfrac{\pi}{8}x^2-\nnfrac{3\pi}{128}x^4-\nnfrac{5\pi}{512}x^6       & T_6^{\mbox{\scriptsize \boldmath $E$},\,0}(x)+\left(1-\nnfrac{175\pi}{512}\right)x^7                           \\[1.75 ex]
 8 & \nnfrac{\pi}{2}-\nnfrac{\pi}{8}x^2-\nnfrac{3\pi}{128}x^4-\nnfrac{5\pi}{512}x^6
   \!-\! \nnfrac{175\pi}{32768}x^8                                                        & T_7^{\mbox{\scriptsize \boldmath $E$},\,0}(x)+\left(1-\nnfrac{175\pi}{512}\right)x^8                           \\[1.75 ex]
 9 & \nnfrac{\pi}{2}-\nnfrac{\pi}{8}x^2-\nnfrac{3\pi}{128}x^4-\nnfrac{5\pi}{512}x^6
   \!-\! \nnfrac{175\pi}{32768}x^8                                                        & T_8^{\mbox{\scriptsize \boldmath $E$},\,0}(x)+\left(1-\nnfrac{11025\pi}{32768}\right)x^9                       \\[1.75 ex]
10 & \nnfrac{\pi}{2}-\nnfrac{\pi}{8}x^2-\nnfrac{3\pi}{128}x^4-\nnfrac{5\pi}{512}x^6
   \!-\! \nnfrac{175\pi}{32768}x^8
   \!-\! \nnfrac{441\pi}{131072}x^{10} \qquad                                             & T_9^{\mbox{\scriptsize \boldmath $E$},\,0}(x)+\left(1-\nnfrac{11025\pi}{32768}\right)x^{10}                    \\[1.75 ex] \hline
\end{array}
\end{equation*}

\vspace*{0.00 mm}

\centerline{\hspace*{-15.0 mm} Table 2.}
%\end{landscape}

\bigskip
\noindent
{\bf \boldmath 3.2 Taylor approximations of $D$-elliptic integrals} %
% %%%%%%%%%%%%%%%%%%%%%%%%%%%%%%%%%%%%%%%%%%%%%%%%%%%%%%%%%%%%%%%%%%%%%%%%%%%%%%%%%%%%%%%%%%%%%%%%%%%%%%%%%%%%%%%%%%%%%%%%%%%%%%%%%%%%%%%%%%%%%%%%%%%%%

\bigskip
%\noindent
For the process of determining a formula for the area $\mathcal{A}_{\mbox{\scriptsize \boldmath $q$}, \,egg}$ of curve
$\mathcal{F}_{\mbox{\scriptsize \boldmath $q$},\, egg}$, Taylor approximations of $D$-elliptic integrals are of
special interest. $D$-elliptic integrals are determined as follows
\begin{equation}
\label{D_def}
\mbox{\boldmath $D$}\left(x\right)
=
\left\{
\begin{array}{ccc}
\dfrac{\mbox{\boldmath $K$}(x) - \mbox{\boldmath $E$}(x)}{x^2} &:& x \in (-1,1) \backslash \{0\} \\[1.5 ex]
0                        &:& x = 0
\end{array}
\right.,
\end{equation}
see formula 8.112/5 in \cite{Ryzhik 2015}.
For $D$-elliptic integrals, the series expansion
\begin{equation}
\label{D_series}
\mbox{\boldmath $D$}\left(x\right)
=
\pi \sum_{i=0}^{\infty}{  \frac{i+1}{2i+1}}\, \left( \frac{(2i+1)!!}{(2i+2)!!} \right)^2 \!\!  x^{2i},
\end{equation}
holds for $x \!\in\! (-1,1)$, see formula 8.115 in \cite{Ryzhik 2015}. Let us note that the series expansion (\ref{D_series})
has a radius of $1$ and that $x = 1$ can not be included in the convergence domain.

\medskip
\noindent
For the elliptic integral $\mbox{\boldmath $D$}(x)$, based on the series expansion (\ref{D_series}), the first and second Taylor approximations
are obtained. The first Taylor approximation is:
\begin{equation}
T_n^{\mbox{\scriptsize \boldmath $D$},\,0}(x)
=
\pi \sum_{i=0}^{[n/2]}{  \frac{i+1}{2i+1}}\, \left( \frac{(2i+1)!!}{(2i+2)!!} \right)^2 \!\!  x^{2i},
\end{equation}
where $x \!\in\! (-1,1)$. For a fixed $\beta \!\in\! (0,1)$, the second Taylor expansion in the right neighborhood of point
$\alpha = 0$ is:
\begin{equation}
\mathbbmsl{T}_n^{\mbox{\scriptsize \boldmath $D$};\,0,\beta}(x) =
    \left\{
    \begin{array}{lll}
    T_{n-1}^{\mbox{\scriptsize \boldmath $D$},\,0}(x) + \dfrac{x^n}{\beta^n}\left(\mbox{\boldmath $D$}(\beta) - T_{n-1}^{\mbox{\scriptsize \boldmath $D$},\,0}(\beta) \right) & : & n \geq 1 \\[1.0 ex]
    \mbox{\boldmath $D$}(\beta) & : & n = 0
    \end{array}
    \right.,
\end{equation}
where $x \!\in\! [-\beta,\beta]$.

\medskip
\noindent
According to \cite{MRL 2019}, the following inequalities hold for the elliptic integral $\mbox{\boldmath $D$}(x):$
\begin{equation}
\label{Nejednakosti_za_D}
\begin{array}{c}
\dfrac{\pi}{4}
\!=\!
T_{0}^{\mbox{\scriptsize \boldmath $D$},\,0}(x) \!=\! T_{1}^{\mbox{\scriptsize \boldmath $D$},\,0}(x) \leq
T_{2}^{\mbox{\scriptsize \boldmath $D$},\,0}(x) \!=\! T_{3}^{\mbox{\scriptsize \boldmath $D$},\,0}(x) \leq
\ldots \leq
T_{2i}^{\mbox{\scriptsize \boldmath $D$},\,0}(x) \!=\! T_{2i+1}^{\mbox{\scriptsize \boldmath $D$},\,0}(x) \leq
\ldots                                                                          \\[1.5 ex]
\leq \mbox{\boldmath $D$}(x) \leq                                                                        \\[1.5 ex]
\ldots \leq
\mathbbmsl{T}_j^{\mbox{\scriptsize \boldmath $D$};\,0,\,\beta}(x) \leq
\ldots               \leq
\mathbbmsl{T}_3^{\mbox{\scriptsize \boldmath $D$};\,0,\,\beta}(x) \leq
\mathbbmsl{T}_2^{\mbox{\scriptsize \boldmath $D$};\,0,\,\beta}(x) \leq
\mathbbmsl{T}_1^{\mbox{\scriptsize \boldmath $D$};\,0,\,\beta}(x) \leq
\mathbbmsl{T}_0^{\mbox{\scriptsize \boldmath $D$};\,0,\,\beta}(x)
\!=\!
\mbox{\boldmath $D$}(\beta)
\end{array}
\end{equation}
for a fixed $\beta \!\in\! (0,1)$ and an arbitrary $x \!\in\! [-\beta,\beta]$.
In the following table, we list the explicit forms of the polynomials $T_j^{{\mbox{\scriptsize \boldmath $D$}}^{\vphantom{2^2}},\,0}(x)$ and $\mbox{$\mathbbmsl{T}$}_j^{\mbox{\scriptsize \boldmath $D$};\,0,\,\beta}(x)$,
for $\beta \!\in\! (0,1)$ and $j = 0, 1, \ldots, 10$:

%\begin{landscape}
%\medskip
\begin{equation*}
\hspace*{-3.5 mm}
\begin{array}{|c|l|l|} \hline
    &
    &                                                                           \\[0.05 ex]
 j & T_j^{{\mbox{\scriptsize \boldmath $D$}}^{\vphantom{2^2}},\,0}(x)                                            & \mbox{$\mathbbmsl{T}$}_j^{\mbox{\scriptsize \boldmath $D$};\,0,\,\beta}(x) \;\mbox{and}\; C_{\mbox{\tiny \boldmath $D$}} = \mbox{\boldmath $D$}(\beta) \;\mbox{for}\; \beta \!\in\! (0,1)                                                       \\[1.5 ex]
   \hline \hline
 0 & \nnfrac{{\pi}^{\vphantom{2}}}{4}
   & C_{\mbox{\tiny \boldmath $D$}}                                                                                 \\[1.5 ex]
 1 & \nnfrac{\pi}{4}
   & T_0^{\mbox{\scriptsize \boldmath $D$};\,0}(x)\!+\!\left(\nnfrac{C_{\mbox{\tiny \boldmath $D$}}}{\beta}\!-\!\nnfrac{\pi}{4\beta}\right)x                              \\[1.5 ex]
 2 & \nnfrac{\pi}{4}\!+\!\nnfrac{3\pi}{32}x^2
   & T_1^{\mbox{\scriptsize \boldmath $D$};\,0}(x)\!+\!\left(\nnfrac{C_{\mbox{\tiny \boldmath $D$}}}{\beta^2}\!-\!\nnfrac{\pi}{4\beta^2}\right)x^2                        \\[1.5 ex]
 3 & \nnfrac{\pi}{4}\!+\!\nnfrac{3\pi}{32}x^2
   & T_2^{\mbox{\scriptsize \boldmath $D$};\,0}(x)\!+\!\left(\nnfrac{C_{\mbox{\tiny \boldmath $D$}}}{\beta^3}\!-\!\nnfrac{\pi}{4\beta^3}\!-\!\nnfrac{3\pi}{32\beta}\right)x^3       \\[1.5 ex]
 4 & \nnfrac{\pi}{4}\!+\!\nnfrac{3\pi}{32}x^2\!+\!\nnfrac{15\pi}{256}x^4
   & T_3^{\mbox{\scriptsize \boldmath $D$};\,0}(x)\!+\!\left(\nnfrac{C_{\mbox{\tiny \boldmath $D$}}}{\beta^4}\!-\!\nnfrac{\pi}{4\beta^4}\!-\!\nnfrac{3\pi}{32\beta^2}\right)x^4     \\[1.5 ex]
 5 & \nnfrac{\pi}{4}\!+\!\nnfrac{3\pi}{32}x^2\!+\!\nnfrac{15\pi}{256}x^4
   & T_4^{\mbox{\scriptsize \boldmath $D$};\,0}(x)\!+\!\left(\nnfrac{C_{\mbox{\tiny \boldmath $D$}}}{\beta^5}\!-\!\nnfrac{\pi}{4\beta^5}\!-\!\nnfrac{3\pi}{32\beta^3}\!-\!\nnfrac{15\pi}{256\beta}\right)x^5   \\[1.5 ex]
 6 & \nnfrac{\pi}{4}\!+\!\nnfrac{3\pi}{32}x^2\!+\!\nnfrac{15\pi}{256}x^4\!+\!\nnfrac{175\pi}{4096}x^6
   & T_5^{\mbox{\scriptsize \boldmath $D$};\,0}(x)\!+\!\left(\nnfrac{C_{\mbox{\tiny \boldmath $D$}}}{\beta^6}\!-\!\nnfrac{\pi}{4\beta^6}\!-\!\nnfrac{3\pi}{32\beta^4}\!-\!\nnfrac{15\pi}{256\beta^2}\right)x^6  \\[1.5 ex]
 7 & \nnfrac{\pi}{4}\!+\!\nnfrac{3\pi}{32}x^2\!+\!\nnfrac{15\pi}{256}x^4\!+\!\nnfrac{175\pi}{4096}x^6
   & T_6^{\mbox{\scriptsize \boldmath $D$};\,0}(x)\!+\!\left(\nnfrac{C_{\mbox{\tiny \boldmath $D$}}}{\beta^7}\!-\!\nnfrac{\pi}{4\beta^7}\!-\!\nnfrac{3\pi}{32\beta^5}\!-\!\nnfrac{15\pi}{256\beta^3}\!-\!\nnfrac{175\pi}{4096\beta}\right)x^7 \\[1.5 ex]
 8 & \nnfrac{\pi}{4}\!+\!\nnfrac{3\pi}{32}x^2\!+\!\nnfrac{15\pi}{256}x^4\!+\!\nnfrac{175\pi}{4096}x^6
   \!\!+\!\! \nnfrac{2205\pi}{65536}x^8
   & T_7^{\mbox{\scriptsize \boldmath $D$};\,0}(x)\!+\!\left(\nnfrac{C_{\mbox{\tiny \boldmath $D$}}}{\beta^8}\!-\!\nnfrac{\pi}{4\beta^8}\!-\!\nnfrac{3\pi}{32\beta^6}\!-\!\nnfrac{15\pi}{256\beta^4}\!-\!\nnfrac{175\pi}{4096\beta^2}\right)x^8 \\[1.5 ex]
 9 & \nnfrac{\pi}{4}\!+\!\nnfrac{3\pi}{32}x^2\!+\!\nnfrac{15\pi}{256}x^4\!+\!\nnfrac{175\pi}{4096}x^6
   \!\!+\!\! \nnfrac{2205\pi}{65536}x^8
   & T_8^{\mbox{\scriptsize \boldmath $D$};\,0}(x)\!+\!\left(\nnfrac{C_{\mbox{\tiny \boldmath $D$}}}{\beta^9}\!-\!\nnfrac{\pi}{4\beta^9}\!-\!\nnfrac{3\pi}{32\beta^7}\!-\!\nnfrac{15\pi}{256\beta^5}\!-\!\nnfrac{175\pi}{4096\beta^3}\!-\!\nnfrac{2205\pi}{65536\beta}\right)x^9  \\[1.5 ex]
10 & \nnfrac{\pi}{4}\!+\!\nnfrac{3\pi}{32}x^2\!+\!\nnfrac{15\pi}{256}x^4\!+\!\nnfrac{175\pi}{4096}x^6
   \!\!+\!\! \nnfrac{2205\pi}{65536}x^8
   \!\!+\!\! \nnfrac{14553\pi}{524288}x^{10} \,\,
   & T_9^{\mbox{\scriptsize \boldmath $D$};\,0}(x)\!+\!\left(\nnfrac{C_{\mbox{\tiny \boldmath $D$}}}{\beta^{10}}\!-\!\nnfrac{\pi}{4\beta^{10}}\!-\!\nnfrac{3\pi}{32\beta^8}\!-\!\nnfrac{15\pi}{256\beta^6}\!-\!\nnfrac{175\pi}{4096\beta^4}\!-\!\nnfrac{2205\pi}{65536\beta^2}\right)x^{10}  \\[1.5 ex]
\hline
\end{array}
\end{equation*}
\centerline{\hspace*{-15.0 mm} Table 3.}
%\end{landscape}

\bigskip

Let us note that there is recent research into elliptic integrals of the first and second kind (and their convexity, monotonicity, approximations, inequalities, applications, ...), see the papers:
\cite{Elliptic 1}, \cite{Elliptic 2}, \cite{Elliptic 3}, \cite{Elliptic 4}, \cite{W-DJ 2023},
\cite{Elliptic 5}, \cite{Elliptic 6}, \cite{Elliptic 7}, \cite{Elliptic 8}, \cite{Elliptic 10}, \cite{AADM 2020},
\cite{Elliptic 11}-\!
%, \cite{Elliptic 12}, \cite{AADM 2019}, \cite{Elliptic 13}, \cite{Elliptic 14} and
\cite{Elliptic 15}.
An approach to the approximate computation of complete elliptic integrals for practical use in water engineering is given in the paper \cite{Elliptic 9}.

\section{Taylor approximations for the area of H\" ugelsch\" affer curve $\mathcal{F}_{\mbox{\scriptsize \boldmath $q$}, \, egg}$}

\bigskip
%\noindent
Let $\mbox{\boldmath $\mathcal{A}$}(k)$ be the size of the area
$\mathcal{A}_{\mbox{\scriptsize \boldmath $q$}, \,egg}$ of curve
$\mathcal{F}_{\mbox{\scriptsize \boldmath $q$}, \, egg}$
as a function of $k \in (0,1)$. For $k =0$, a degenerative case is obtained
wherein the egg-shaped part of the curve is reduced to an ellipse. Then, we define that $\mbox{\boldmath $\mathcal{A}$}(0) = a b \mbox{\boldmath $q$} \pi$.
Additionally, for $k = 1$, another degenerative case is obtained wherein the egg-shaped part of the curve
is bounded by a part of a parabola and a line. Then, we define that $\mbox{\boldmath $\mathcal{A}$}(1) = \frac{8}{3} a b \mbox{\boldmath $q$}$.
Based on the degenerative cases and the formula for the non-degenerative case (\ref{A1&A2}), with the use of the
{\boldmath $D$}-elliptic integral, it follows that
\begin{equation}
\label{3_putica}
\mbox{\boldmath $\mathcal{A}$}(k)
=
\left\{
\begin{array}{ccc}
a b \mbox{\boldmath $q$} \pi &:& k = 0, \\[1.0 ex]
\dfrac{4}{3}\, a \, b \,\mbox{\boldmath $q$}\,
{\bigg (}\mbox{\boldmath $K$}(k)
+
\mbox{\boldmath $E$}(k)
-
\mbox{\boldmath $D$}\left(k\right) {\bigg )} &:& k \in (0,1), \\[1.0 ex]
\dfrac{8}{3} a b \mbox{\boldmath $q$} &:& k = 1.
\end{array}
\right.
\end{equation}
For $k \in (0,1)$, using the formulae (\ref{E_series}), (\ref{K_series}) i (\ref{D_series}) we obtain a series expansion:
\begin{equation}
\label{A_series}
\mbox{\boldmath $\mathcal{A}$}(k)
=
a \, b \,\mbox{\boldmath $q$}\, \pi
-
a \, b \,\mbox{\boldmath $q$}\, \pi
\sum_{i=1}^{\infty}{ \frac{1}{(2i-1)(i+1)}}\left( \frac{(2i-1)!!}{(2i)!!} \right)^2 \!\!  k^{2i}.
\end{equation}
Let us note that the series expansion (\ref{A_series}) has a radius of $1$ and that the value $k = 1$ can be included in the
convergence domain in accordance with Raabe's test. Furthermore, for $k = 0$ and for $k = 1$ the formulae (\ref{3_putica}) and
(\ref{A_series}) give the same results.

\bigskip
\begin{remark}
Based on the two previous expressions, (\ref{3_putica}) and (\ref{A_series}), for $k = 1$:
\begin{equation}
\label{EPi1}
\mbox{\boldmath $\mathcal{A}$}(1)
=
\dfrac{8}{3} a b \mbox{\boldmath $q$}
\end{equation}
and
\begin{equation}
\label{EPi2}
\mbox{\boldmath $\mathcal{A}$}(1)
=
a \, b \,\mbox{\boldmath $q$}\, \pi
-
a \, b \,\mbox{\boldmath $q$}\, \pi
\sum_{i=1}^{\infty}{ \frac{1}{(2i-1)(i+1)}}\left( \frac{(2i-1)!!}{(2i)!!} \right)^2
\end{equation}
the following representation holds:
\begin{equation}
\label{EPi}
\dfrac{1}{\pi}
=
\dfrac{3}{8}\left(1-\sum_{i=1}^{\infty}{ \frac{1}{(2i-1)(i+1)}}\left( \frac{(2i-1)!!}{(2i)!!} \right)^2\right).
\end{equation}
\end{remark}

\medskip
\noindent
The approximation of the numbers $\pi$ and $\frac{1}{\pi}$ was the topic of research in the papers \cite{Pi 1}, \cite{Pi 2},  \cite{Pi 3}, \cite{Pi 4}, \cite{Pi 6} and \cite{Pi 5}.

\medskip
\noindent
For $\mbox{\boldmath $\mathcal{A}$}(k)$, based on the series expansion (\ref{A_series}), the first and second Taylor approximations are obtained, which we denote with the formulae (\ref{FirstA}) and (\ref{SecondA}) in the remainder of this paper.

\smallskip
\noindent
The first Taylor approximation is:
\begin{equation}
\label{FirstA}
T_n^{\mbox{\scriptsize \boldmath $\mathcal{A}$},\,0}(k)
=
a \, b \,\mbox{\boldmath $q$}\, \pi
-
a \, b \,\mbox{\boldmath $q$}\, \pi
\sum_{i=1}^{[n/2]}{ \frac{1}{(2i-1)(i+1)}}\, \left( \frac{(2i-1)!!}{(2i)!!} \right)^2 \!\! k^{2i}, \,\,\, k \!\in\! (0,1).
\end{equation}
For a fixed $\beta \!\in\! (0,1]$, the second Taylor approximation in the right neighborhood of point
$\alpha = 0$ is:
\begin{equation}
\label{SecondA}
\mathbbmsl{T}_n^{\mbox{\scriptsize \boldmath $\mathcal{A}$};\,0,\beta}(k) =
    \left\{
    \begin{array}{lll}
    T_{n-1}^{\mbox{\scriptsize \boldmath $\mathcal{A}$},\,0}(k) + \dfrac{x^n}{\beta^n}\left(\mbox{\boldmath $\mathcal{A}$}(\beta) - T_{n-1}^{\mbox{\scriptsize \boldmath $\mathcal{A}$},\,0}(\beta) \right) & : & n \geq 1 \\[1.0 ex]
    \mbox{\boldmath $\mathcal{A}$}(\beta) & : & n = 0
    \end{array}
    \right.\!\!, \,\,\,k \!\in\! (0,\beta].
\end{equation}

\bigskip
\noindent
According to \cite{MRL 2019}, the following inequalities hold:
\begin{equation}
\label{Nejednakosti_za_A}
\begin{array}{c}
a \, b \,\mbox{\boldmath $q$}\,{\pi}^{\vphantom{2^2}}
\!=\!
T_{0}^{\mbox{\scriptsize \boldmath $\mathcal{A}$},\,0}(k) \!=\! T_{1}^{\mbox{\scriptsize \boldmath $\mathcal{A}$},\,0}(k) \geq
T_{2}^{\mbox{\scriptsize \boldmath $\mathcal{A}$},\,0}(k) \!=\! T_{3}^{\mbox{\scriptsize \boldmath $\mathcal{A}$},\,0}(k) \geq
\ldots \geq
T_{2i}^{\mbox{\scriptsize \boldmath $\mathcal{A}$},\,0}(k) \!=\! T_{2i+1}^{\mbox{\scriptsize \boldmath $\mathcal{A}$},\,0}(k) \geq
\ldots                                                                          \\[1.5 ex]
\geq \mbox{\boldmath $\mathcal{A}$}(k) \geq                                                                        \\[1.5 ex]
\ldots \geq
\mathbbmsl{T}_j^{\mbox{\scriptsize \boldmath $\mathcal{A}$};\,0,\,\beta}(k) \geq
\ldots               \geq
\mathbbmsl{T}_3^{\mbox{\scriptsize \boldmath $\mathcal{A}$};\,0,\,\beta}(k) \geq
\mathbbmsl{T}_2^{\mbox{\scriptsize \boldmath $\mathcal{A}$};\,0,\,\beta}(k) \geq
\mathbbmsl{T}_1^{\mbox{\scriptsize \boldmath $\mathcal{A}$};\,0,\,\beta}(k) \geq
\mathbbmsl{T}_0^{\mbox{\scriptsize \boldmath $\mathcal{A}$};\,0,\,\beta}(k)
\!=\!
\mbox{\boldmath $\mathcal{A}$}(\beta)
\end{array}
\end{equation}
for a fixed $\beta \!\in\! (0,1)$ and an arbitrary $k \!\in\! (0, \beta]$.

\medskip
\noindent
We specify a list of {\sc Taylor} approximations for $\mathcal{A}(k)$ where $k = \dfrac{\mbox{\boldmath $q$}^2 w}{a}$, i.e. $k \!\in\! (0, 1]$, in table 4:
%\begin{landscape}
\medskip
\begin{equation*}
\begin{array}{|c|l|l|} \hline
    &
    &                                                                           \\[0.05 ex]
 j & T_j^{{\mbox{\scriptsize \boldmath $\mathcal{A}$}}^{\vphantom{2^2}},\,0}(k)
   & \mbox{$\mathbbmsl{T}$}_j^{\mbox{\scriptsize \boldmath $\mathcal{A}$};\,0,\,\beta}(k) \;\mbox{and}\; \beta =1                                          \\[1.5 ex]
   \hline \hline
    &
    &                                                                           \\[0.05 ex]
 0  & a \, b \,\mbox{\boldmath $q$}\,{\pi}^{\vphantom{2^2}}
    & \nnfrac{8}{3}a \, b \,\mbox{\boldmath $q$}                                                  \\[1.5 ex]
 1  & a \, b \,\mbox{\boldmath $q$}\,\pi
    & T_0^{\mbox{\scriptsize \boldmath $\mathcal{A}$},\,0}(k)+ a \, b \,\mbox{\boldmath $q$}\,\left(\nnfrac{8}{3}-\pi\right)k                    \\[1.5 ex]
 2  & a \, b \,\mbox{\boldmath $q$}\,\pi\, \left(1-\nnfrac{1}{8}k^2\right)
    & T_1^{\mbox{\scriptsize \boldmath $\mathcal{A}$},\,0}(k)+ a \, b \,\mbox{\boldmath $q$}\,\left(\nnfrac{8}{3}-\pi\right)k^2                  \\[1.5 ex]
 3  & a \, b \,\mbox{\boldmath $q$}\,\pi\, \left(1-\nnfrac{1}{8}k^2\right)
    & T_2^{\mbox{\scriptsize \boldmath $\mathcal{A}$},\,0}(k)+ a \, b \,\mbox{\boldmath $q$}\,\left(\nnfrac{8}{3}-\nnfrac{7}{8}\pi\right)k^3     \\[1.5 ex]
 4  & a \, b \,\mbox{\boldmath $q$}\,\pi\, \left(1-\nnfrac{1}{8}k^2-\nnfrac{1}{64}k^4\right)
    & T_3^{\mbox{\scriptsize \boldmath $\mathcal{A}$},\,0}(k)+ a \, b \,\mbox{\boldmath $q$}\,\left(\nnfrac{8}{3}-\nnfrac{7}{8}\pi\right)k^4     \\[1.5 ex]
 5  & a \, b \,\mbox{\boldmath $q$}\,\pi\, \left(1-\nnfrac{1}{8}k^2-\nnfrac{1}{64}k^4\right)
    &T_4^{\mbox{\scriptsize \boldmath $\mathcal{A}$},\,0}(k)+ a \, b \,\mbox{\boldmath $q$}\,\left(\nnfrac{8}{3}-\nnfrac{55}{64}\pi\right)k^5    \\[1.5 ex]
 6  & a \, b \,\mbox{\boldmath $q$}\,\pi\, \left(1-\nnfrac{1}{8}k^2-\nnfrac{1}{64}k^4-\nnfrac{5}{1024}k^6\right)
    & T_5^{\mbox{\scriptsize \boldmath $\mathcal{A}$},\,0}(k)+ a \, b \,\mbox{\boldmath $q$}\,\left(\nnfrac{8}{3}-\nnfrac{55}{64}\pi\right)k^6   \\[1.5 ex]
 7  & a \, b \,\mbox{\boldmath $q$}\,\pi\, \left(1-\nnfrac{1}{8}k^2-\nnfrac{1}{64}k^4-\nnfrac{5}{1024}k^6\right)
    & T_6^{\mbox{\scriptsize \boldmath  $\mathcal{A}$},\,0}(k)+ a \, b \,\mbox{\boldmath $q$}\,\left(\nnfrac{8}{3}-\nnfrac{875}{1024}\pi\right)k^7 \\[1.5 ex]
 8   & a \, b \,\mbox{\boldmath $q$}\,\pi\, \left(1-\nnfrac{1}{8}k^2-\nnfrac{1}{64}k^4-\nnfrac{5}{1024}k^6
   \!-\! \nnfrac{35}{16384}k^8\right)
    & T_7^{\mbox{\scriptsize \boldmath $\mathcal{A}$},\,0}(k)+ a \, b \,\mbox{\boldmath $q$}\,\left(\nnfrac{8}{3}-\nnfrac{875}{1024}\pi\right)k^8 \\[1.5 ex]
 9 & a \, b \,\mbox{$q$}\,\pi\, \left(1-\nnfrac{1}{8}k^2-\nnfrac{1}{64}k^4-\nnfrac{5}{1024}k^6
   \!-\! \nnfrac{35}{16384}k^8\right)
   & T_8^{\mbox{\scriptsize \boldmath $\mathcal{A}$},\,0}(k)+ a \, b \,\mbox{$q$}\,\left(\nnfrac{8}{3}-\nnfrac{13965}{16384}\pi\right)k^9  \\[1.5 ex]
10 & a \, b \,\mbox{$q$}\,\pi\, \left(1-\nnfrac{1}{8}k^2-\nnfrac{1}{64}k^4-\nnfrac{5}{1024}k^6
   \!-\! \nnfrac{35}{16384}k^8 \!-\! \nnfrac{147}{131072}k^{10}\right)\,
   & T_9^{\mbox{\scriptsize \boldmath $\mathcal{A}$},\,0}(k)+ a \, b \,\mbox{$q$}\,\left(\nnfrac{8}{3}-\nnfrac{13965}{16384}\pi\right)k^{10}  \\[1.5 ex]
\hline
\end{array}
\end{equation*}

\centerline{\hspace*{-15.0 mm} Table 4.}
%\end{landscape}

\section{\boldmath Approximative formulae for the area of curve $\mathcal{F}_{\mbox{\scriptsize \boldmath $q$}, \, egg}$}

\medskip
\noindent
\begin{theorem} \label{Theorem P-Egg-E}
(Theorem for the area of H\" ugelsch\" affer egg curves).
The following two estimations hold:
\begin{equation}
\label{procena1}
\nnfrac{8}{3}a \, b \,\mbox{\boldmath $q$}
\leq
\mathcal{A}(k)
\leq
\pi \, a \, b \,\mbox{\boldmath $q$}
\end{equation}
and
\begin{equation}
\label{procena2}
\nnfrac{8}{3}a \, b \,\mbox{\boldmath $q$} + \Delta_{\mbox{\scriptsize \boldmath $q$}}
\leq
\mathcal{A}(k)
\leq
\pi \, a \, b \,\mbox{\boldmath $q$}-\nabla_{\mbox{\scriptsize \boldmath $q$}}
\end{equation}
where
\begin{equation}
\Delta_{\mbox{\scriptsize \boldmath $q$}} =
a \, b \,\mbox{\boldmath $q$}\,\pi \, (1-k)=
\left\{
\begin{array}{ccc}
b \,\pi \, (a-w)                        &:& w < a \\[0.8 ex]
\nfrac{a^2} {w^2} \, b \,\pi \, (w-a)   &:& w > a
\end{array}
\right. ,
\qquad
\nabla_{\mbox{\scriptsize \boldmath $q$}} =
\nfrac{\pi}{8}\,a \, b \,\mbox{\boldmath $q$}\,k^2 =
\left\{
\begin{array}{ccc}
\nfrac{\pi} {8\,a} \, b  \, w^2     &:& w < a \\[0.8 ex]
\nfrac{\pi} {8\,w^3} \, a^4 \, b    &:& w > a
\end{array}
\right. \!
\end{equation}
for $k \in [0,1]$.
\end{theorem}

\noindent
\textit{Proof}.
The series expansion $\mathcal{A}(k)$ given by (\ref{A_series}) is a monotonously decreasing function when $k \in [0,1]$ (based on the term by term differentiation of the series).
Thus, the estimate (\ref{procena1}) is obtained based on (\ref{3_putica}) and (\ref{A_series}).

\bigskip
\noindent
Based on (\ref{Nejednakosti_za_A}) and table 4, the following is a proof of the estimate (\ref{procena2}):
\begin{equation}
\nnfrac{8}{3}a \, b \,\mbox{\boldmath $q$}
\leq
a \, b \,\mbox{\boldmath $q$}\, \pi + a \, b \,\mbox{\boldmath $q$}\,\left( \nnfrac{8}{3}-\pi \right)k
\leq
\mathcal{A}(k)
\leq
\pi \, a \, b \,\mbox{\boldmath $q$}\,\left( 1-\nnfrac{1}{8}k^2\right)
\leq
\pi \, a \, b \,\mbox{\boldmath $q$}
\end{equation}
i.e.
\begin{equation}
\nnfrac{8}{3}a \, b \,\mbox{\boldmath $q$}
\leq
 \nnfrac{8}{3}a \, b \,\mbox{\boldmath $q$}\, + \, a \, b \,\mbox{\boldmath $q$}\,\pi \, (1-k)
\leq
\mathcal{A}(k)
\leq
\pi \, a \, b \,\mbox{\boldmath $q$} -\nnfrac{\pi}{8} \, a \, b \,\mbox{\boldmath $q$}\,k^2
\leq
\pi \, a \, b \,\mbox{\boldmath $q$}.
\end{equation}

\begin{figure}[htb]
\vspace{-0.3cc}
  \centerline{
    \includegraphics*[width=0.6\textwidth]{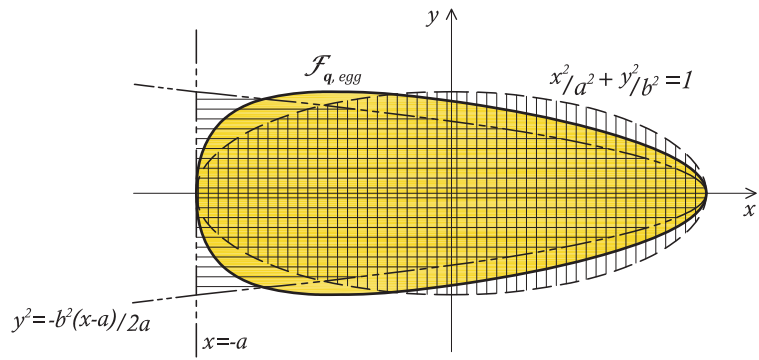}
  }
\caption{A comparison of the areas of a parabola, H\" ugelsch\" affer   egg curve $\mathcal{F}_{\mbox{\scriptsize \boldmath $q$},\,egg}$, and an ellipse; \, Source: \copyright\,First Author}
\label{Fig.3}       % Give a unique label
\end{figure}

\noindent
The estimate (\ref{procena1}) is graphically illustrated in Fig. \ref{Fig.3}.

%% FIGURE 3 %%%%%%%%%%%%%%%%%%%%%%%%%%%%%%%%%%%%%%%%%%%%%%%%%%%%%%%%%%%%%%%%%%%%%%%%%%%%%%%%%%%%%%%%

\begin{remark}
Based on the series of inequalities (\ref{Nejednakosti_za_A}), it is possible to obtain even better estimates for $\mathcal{A}(k)$ with the appropriate polynomials given in table 4.
\end{remark}

\section{Conclusion}

Applying Taylor series and double Taylor series (Section 3), we have obtained novel approximations for
$\mbox{\boldmath $K$}$, $\mbox{\boldmath $E$}$ and $\mbox{\boldmath $D$}$ elliptic integrals (see table 1 -- 3),
as well as approximative formulae for calculating the area $\mathcal{A}(k)$ of the egg-shaped part of H\" ugelsch\" affer curves (see table 4 in Section 4).
Furthermore, based on the expression (\ref{EPi1}) which represents the area of a part of a parabola (see Fig. \ref{Fig.3}),
and the expression (\ref{EPi2}) which represents the series expansion $\mathcal{A}(k)$ when $k=1$,
a new representation (\ref{EPi}) of the number $1/\pi$ has arisen.

With the development of new software tools and applets (see \cite{Applet 2022}), the use of the newly-introduced formulae
%(see Fig. \ref{Fig.3}, golden colored area)
(\ref{A1&A2}) or (\ref{3_putica}) for the calculation of the area of the egg-shaped part of H\" ugelsch\" affer curves $\mathcal{A}_{\mbox{\scriptsize \boldmath $q$},\,egg}$ or
% (aero-, civil-, water-, bio-, traffic-, food-, agriculture-, ...),\,
$\mathcal{A}(k)$ could be significant in various areas of engineering, poultry industry and ornithology.

\end{document}